\documentclass[12pt]{article}
\usepackage{amssymb}
\usepackage{amsmath, amscd}
\usepackage{color}
\usepackage{latexsym}
\usepackage[dvips]{graphicx}
\usepackage{graphicx}
\usepackage{hyperref}
\usepackage{mathrsfs}
\usepackage{changes}
\newtheorem{theorem}{Theorem}[section]
\newtheorem{lemma}[theorem]{Lemma}
\newtheorem{corollary}[theorem]{Corollary}
\newtheorem{definition}[theorem]{Definition}
\newtheorem{notations}[theorem]{Notations}
\newtheorem{algorithm}[theorem]{Algorithm}

\newcommand{\qed}{\hfill $\Box$ }

\newcommand{\proof}{\noindent{\bf Proof.}\ \ }

\usepackage{tikz} 
\usetikzlibrary{arrows,automata,positioning}
\usetikzlibrary{graphs,graphs.standard,math}
\usetikzlibrary{shapes.geometric}
\usepackage{textpos}
\usetikzlibrary {quotes}
\usepackage{verbatim}
\usetikzlibrary{arrows,shapes}
\begin{document}

\title{\Large {\bf A decomposition structure of resonance graphs that are daisy cubes}}
\maketitle

\begin{center}
{\large \bf Zhongyuan Che$^{a}$, Zhibo Chen$^b$} 
\end{center}


\smallskip 

\begin{center}{$^a$ \it Department of Mathematics, Penn State University, Beaver Campus, Monaca, PA 15061,  zxc10@psu.edu} \end{center}

\begin{center}{$^b$ \it Department of Mathematics, Penn State University, Greater Allegheny Campus, McKeesport, PA 15132,  zxc4@psu.edu}\end{center}

\smallskip 
\begin{center}{\large Dedicated to Professor Fuji Zhang on the occasion of his 88th birthday}\end{center}

\begin{abstract}

It has recently been shown in [\emph{Discrete Appl. Math.} {\bf 366} (2025) 75--85]  that the resonance graph  
of a plane elementary bipartite graph $G$  is a daisy cube if and only if $G$ is  peripherally 2-colorable.
Let $G$ be a peripherally 2-colorable graph and $R(G)$ be its resonance graph.
We provide a decomposition structure of $R(G)$ with respect to an arbitrary finite face of $G$
together with a proper labelling for the vertex set of $R(G)$.
An algorithm is obtained to generate a proper labelling for all perfect matchings of $G$  
which induces an  isometric embedding of $R(G)$ as a daisy cube into an $n$-dimensional hypercube,
where $n$ is the isometric dimension of $R(G)$.
Moreover, the algorithm can be applied to generate such a proper labelling for all perfect matchings 
of any plane weakly elementary bipartite graph whose each elementary component with more than two vertices is peripherally 2-colorable.
We also compare two binary codings for all perfect matchings of $G$ which induces distinct structures on $R(G)$:
one as  a daisy cube and the other as a finite distributive, respectively.

\vskip 0.2in
\noindent {\emph{keywords}}: daisy cube, isometric embedding, peripherally 2-colorable graph,
peripheral convex $\le$-expansion,  plane (weakly) elementary bipartite graph, proper labelling, resonance graph, $Z$-transformation graph

\end{abstract}

\section{Introduction} 
A \textit{perfect matching} (or, a \textit{1-factor}) of a graph is a set of vertex disjoint edges  incident to all vertices of the graph.
All graphs considered in this paper have a perfect matching unless specified otherwise. 
A bipartite graph is \textit{elementary} if and only if it is connected and each edge is contained in a perfect matching of the graph \cite{LP86}. 
\textit{Elementary components} of a bipartite graph 
are components of the subgraph obtained by removing all edges not contained in any perfect matching of the graph.
A plane bipartite graph (not necessarily connected) is called \textit{weakly elementary} 
if deleting all edges not contained in any perfect matching of the graph does not result in any new finite faces \cite{Z06}. 
It is easily seen that any plane elementary bipartite graph is also weakly elementary.

The \textit{resonance graph}  of a plane bipartite graph $G$, denoted by $R(G)$,  is a graph 
whose vertices are perfect matchings of $G$, and two perfect matchings $M_1,M_2$ are adjacent  in $R(G)$ 
if their symmetric difference $M_1 \oplus M_2$ (that is, the set of edges contained in either $M_1$ or $M_2$ but not both)
forms exactly one cycle that is the  periphery of a finite face $s$ of $G$, 
and we say that edge $M_1M_2$ of $R(G)$ has the \textit{face-label} $s$.
The concept of resonance graphs was first introduced by chemists  \cite{G82, G83, RKEC96}, and it was shown \cite{R97} 
that the leading eigenvalues of the resonance graphs correlate to the resonance energy of benzenoids.
Professor F. Zhang et al.  introduced the concept of resonance graphs in terms of the $Z$-transformation graphs
of hexagonal systems \cite{ZGC88}, and extended the concept to the $Z$-transformation graphs of plane bipartite graphs \cite{ZZ00}.
One fundamental result \cite{F03, ZZY04} of the concept is that the resonance graph of a plane bipartite graph $G$ is connected 
if and only if $G$ is weakly elementary.
More results on resonance graphs have been developed since the early survey paper  \cite{Z06} by H. Zhang  in 2006.
Resonance graphs display rich structures which are related to 
finite distributive lattices  \cite{C24+, LZ03, Z10}, median graphs  \cite{ZLS08}, 
Fibonacci cubes \cite{KZ05, ZOY09}, and daisy cubes \cite{BCTZ24+, BCTZ25, Z18}.
Recently,  peripherally 2-colorable graphs were introduced \cite{BCTZ25} as a special type of plane elementary bipartite graphs 
to characterize when a resonance graph is a daisy cube.
It was shown \cite{BCTZ25} that if $G$ is a plane elementary bipartite graph, 
then its resonance graph $R(G)$  is a daisy cube if and only if $G$  is  peripherally 2-colorable.
Furthermore, if $G$ is a plane bipartite graph, then $R(G)$  is a daisy cube if and only if $G$ is weakly elementary 
and each elementary component of $G$ with more than two vertices is  peripherally 2-colorable.

The concept of a reducible face of a plane elementary bipartite graph introduced
by  H. Zhang and F. Zhang \cite{ZZ00} has played a key role in the study of a decomposition structure of a resonance graph. 
Let $G$ be a plane elementary bipartite graph. 
A finite face $s$ of $G$  is called a \textit{reducible face}  if the
common periphery of $s$ and $G$ is an odd length path $P$, and
the removal of the internal vertices and edges of $P$ results in a plane elementary bipartite graph. 
A decomposition structure of the resonance graph $R(G)$ 
with respect to a reducible face of $G$ was provided in \cite{C18}, 
and applied frequently \cite{BCTZ24+, BCTZ25, BCTZ24, C19, C20, C21, C24+} 
to analyze the structural properties of resonance graphs. 

In this paper, we focus on considering peripherally 2-colorable graphs.  
Let $G$ be a peripherally 2-colorable graph. Then its resonance graph $R(G)$ is a daisy cube, 
and there is a 1--1 correspondence 
between the set of finite faces of $G$ and the set of  $\Theta$-classes of $R(G)$ such that
each finite face of $G$ is the face-label of  a unique $\Theta$-class of $R(G)$ \cite{BCTZ25}.
Motivated from the fact that any $\Theta$-class of a daisy cube is peripheral \cite{T20}, 
we first provide a  decomposition structure of $R(G)$ with respect to any finite face of $G$.
A decomposition structure of $R(G)$ with respect to a reducible face of $G$ follows as a corollary.
This allows us to construct $R(G)$  from the one-edge graph by a sequence of peripheral convex $\le$-expansions 
with respect to a reducible face decomposition of  $G$ so that a proper labelling for perfect matchings is updated 
at each expansion step which induces an isometric embedding 
of the newly generated resonance graph as a daisy cube into a hypercube.
Algorithm \ref{A:DaisyCube} is obtained naturally to generate a proper labelling for all perfect matchings of $G$ which induces an 
isometric embedding of $R(G)$  as a daisy cube into an $n$-dimensional hypercube, 
where $n$ is both the isometric dimension of $R(G)$
and the number of finite faces of $G$. 
Furthermore, Algorithm \ref{A:DaisyCube} can be applied to generate a proper labelling for all perfect matchings of a
plane weakly elementary bipartite graph $G'$ whose elementary components with more than two vertices are all peripherally 2-colorable,
and the proper labelling induces an  isometric embedding of $R(G')$ as a daisy cube into a $d$-dimensional hypercube,
where $d$ is the isometric dimension of $R(G')$ as well as the summation of the numbers of finite faces over all peripherally 2-colorable components of $G'$.
We conclude the paper by comparing two distinct binary codings for all perfect matchings of $G$ which 
induce $R(G)$  as  a daisy cube and as a finite distributive, respectively.

\section{Preliminaries}\label{S:Pre}

\subsection{Convex expansions and median graphs}
Let $G$ be a graph with the vertex set $V(G)$ and the edge set $E(G)$. 
An induced subgraph of $G$ generated by  $X \subseteq V(G)$ is denoted as $\langle X \rangle$. 
Let $u, v$ be two vertices of $G$. The distance between $u,v$ in $G$,  denoted by $d_G(u,v)$, 
 is the length of a shortest path between $u$ and $v$. 
 The set of all vertices that are on some shortest path between $u, v$ in $G$ is denoted by $I_G(u,v)$.
Two edges $x_1x_2$ and $y_1y_2$ of a connected graph $G$ are  
in \textit{relation $\Theta$} if $d_G(x_1, y_1) + d_G(x_2, y_2) \neq  d_G(x_1, y_2) + d_G(x_2, y_1)$.
Relation $\Theta$ is reflexive and symmetric but not necessarily transitive. 
Let  $xy$ be an edge of a connected graph $G$. Then four subsets of the vertex set $V(G)$ 
and one subset of the edge set $E(G)$ can be defined as follows.
\begin{eqnarray*}
W_{xy} &=& \{w \mid w \in V(G), d_G(w,x) < d_G(w,y)\},\\
W_{yx} &=& \{w \mid w \in V(G), d_G(w,y) < d_G(w,x)\},\\
U_{xy} &=& \{u \in W_{xy} \mid u \mbox{ is adjacent to a vertex in } W_{yx}\},\\
U_{yx} &=& \{v \in W_{yx} \mid v \mbox{ is adjacent to a vertex in } W_{xy}\},\\
F_{xy}&=&\{ e \in E(G) \mid  e \Theta xy\}.
\end{eqnarray*}

It is obvious that $W_{xy} \cap W_{yx} =\emptyset$, and so $U_{xy} \cap U_{yx} =\emptyset$. 
By Proposition 11.7 in \cite{HIK11}, if $G$ is a connected bipartite graph, then 
the spanning subgraph $G-F_{xy}$ has exactly two components, 
namely, $\langle W_{xy} \rangle$ and  $\langle W_{yx} \rangle$.

 Let $H$ be an induced subgraph of a connected graph $G$. Then $H$ is \textit{isometric} in $G$ if 
$d_G(x,y)=d_H(x,y)$ for any two vertices $x, y$ of $H$; and $H$ is \textit{convex} in $G$ if
all vertices on any shortest path between $x$ and $y$ in $G$ are also in $H$  for any two vertices $x$ and $y$ of $H$. 
 A connected graph $G$ is called a \textit{median graph} if $|I_G(u,v) \cap I_G(v,w) \cap I_G(u,w)|=1$ 
for every three vertices $u, v, w$ of $G$.

Let $n$ be a positive integer. Let $\mathcal{B}^{n}$ be the set of all \textit{binary strings} of length $n$. 
A \textit{hypercube} $Q_n$ of dimension $n$ is a graph with the vertex set  $\mathcal{B}^{n}$
such that two vertices are adjacent if they differ in exactly one position. 
A \textit{partial cube} is  an isometric subgraph of a hypercube.  
 All median graphs are partial cubes.
The \textit{isometric dimension} of a partial cube $G$, denoted by $\mathrm{idim}(G)$,  
is the least integer $n$ for which $G$ embeds isometrically into a hypercube $Q_n$.
For any partial cube $G$, relation $\Theta$ is transitive and forms an equivalence relation on $E(G)$,
and $\mathrm{idim}(G)$ is the number of $\Theta$-classes of $G$ \cite{D73}.
By Theorem 11.8 in \cite{HIK11}, if $G$ is a partial cube, then for any edge $xy$ of $G$, 
$F_{xy}$ is an equivalence $\Theta$-class.
So, any $\Theta$-class of a partial cube $G$ can be represented as $F_{xy}$ for some edge $xy \in E(G)$.
Moreover, if  $xy$ is an edge of a partial cube $G$ 
such that vertex $x \in \mathcal{B}^{n}$ (respectively,  $y \in \mathcal{B}^{n}$) has $0$ (respectively, $1$) at position $i$
for some $i \in \{1, 2, \ldots, n\}$, then 
$W_{xy}$ (respectively, $W_{yx}$) contains all vertices of $G$ in $\mathcal{B}^{n}$ having  $0$  (respectively, $1$) at position $i$.

\begin{theorem} \label{T:Median} \cite{HIK11}
Let $xy$ be an edge of a connected bipartite graph $G$. 
Then $G$ is a median graph if and only if the following three conditions are satisfied:

(i) $F_{xy}$ is a matching defining an isomorphism between $\langle U_{xy} \rangle$ and $\langle U_{yx} \rangle$,

(ii) $\langle U_{xy} \rangle$ is convex in $\langle W_{xy} \rangle$ and $\langle U_{yx} \rangle$ is convex in $\langle W_{yx} \rangle$,

(iii) $\langle W_{xy} \rangle$ and $\langle W_{yx} \rangle$ are median graphs.
\end{theorem}

If $xy$ is an edge of a partial cube $G$ such that $U_{xy}=W_{xy}$, 
then  $\langle W_{xy} \rangle$ is called a \textit{peripheral subgraph} of $G$.
A $\Theta$-class $F_{xy}$ of a partial cube $G$ is called a \textit{peripheral $\Theta$-class} 
if at least one of  $\langle W_{xy} \rangle$ and  $\langle W_{yx} \rangle$ is a peripheral subgraph of $G$.

Let $G_1$ and $G_2$ be isometric subgraphs of $G$ such that $V(G_1) \cup V(G_2)=V(G)$
and $V(G_1) \cap V(G_2) \neq \emptyset$.
Assume that  there are no edges between $V(G_1) \setminus V(G_2)$ and $V(G_2)\setminus V(G_1)$.
Take disjoint copies of $G_1$ and $G_2$,  
and add edges between $G_1$ and $G_2$ such that every vertex of $V(G_1) \cap V(G_2)$ in $G_1$ 
is connected with the same vertex of $V(G_1) \cap V(G_2)$ in $G_2$. 
Then the resulted graph is called an \textit{expansion} of $G$.
If $\langle V(G_1) \cap V(G_2) \rangle$ is a convex  subgraph of $G$, then the expansion is called a \textit{convex expansion}.
If $G_1=G$, then the expansion is called \textit{peripheral expansion} of $G$ with respect to $G_2$, 
and denoted by $\mathrm{pe}(G,G_2)$.

A graph is a partial cube if and only if it can be obtained from the one-vertex graph by a sequence of expansions \cite{C88}. 
A  graph is a median graph if and only if it can be obtained from the one-vertex graph by a sequence of convex expansions \cite{M78}.

\subsection{Daisy cubes}

Let $n$ be a positive integer. Let $\mathbb{B}^n=(\mathcal{B}^n, \le)$ 
be a poset on the set  $\mathcal{B}^{n}$  of all \textit{binary strings} of length $n$
such that $u_1u_2 \ldots u_n \leq v_1v_2 \ldots v_n$ if $u_i \leq v_i$ for all $1 \le i \le n$.
A \textit{daisy cube}  is an induced subgraph of $Q_n$ generated by 
a  subset $X \subset \mathcal{B}^n$, and 
can be represented as $Q_n(X)=\langle \{ u \in \mathcal{B}^{n} \ | \ u \leq x \textrm{ for some } x \in X \} \rangle.$
Daisy cubes are partials cubes \cite{KM19}. 
It was pointed out in \cite{T20} that if a graph $K$ is isomorphic to a daisy cube of isometric dimension $n$, 
then there can be more than one isometric embedding of $K$ into the hypercube $Q_n$.
A binary coding for the vertex set of $K$ is called a \textit{proper labelling}
if it induces an isometric embedding of $K$ as a daisy cube into $Q_n$. 
 
 Three terminologies, namely an operator $o$  (Definition 2.5  \cite{T20}),  a $\le$-subgraph (Definition 2.6 \cite{T20}),
and a $\leq$-expansion, were introduced  for daisy cubes in \cite{T20}. 
We extend these concepts from daisy cubes to partial cubes.

\begin{definition}\label{D:Definitions} 
Let $G$ be a partial cube whose vertex set is contained in $(\mathcal{B}^n, \le)$. 
Let $\mathscr{P} (V (G))$ be the set of subsets of $V(G)$.
Then the \textit{operator} $o : \mathscr{P} (V (G)) \rightarrow \mathscr{P} (V (G))$
is defined by $o(X) =\{u \in V(G) | u  \le v$  for some $v \in X\}$.
 An induced subgraph $H$ of $G$ is called a \textit{$\le$-subgraph} of $G$
if $V (H) = o(V (H))$. Let $H$ be a $\le$-subgraph of $G$. 
 Then a peripheral expansion $\mathrm{pe}(G, H)$ is called a \textit{$\le$-expansion}.
\end{definition}
  
It was shown in \cite{T20} that a connected graph  is a daisy cube if and only if 
it can be obtained from the one-vertex graph by a sequence of $\le$-expansions.
We organize the decomposition structure of a daisy cube together with a proper labelling for its vertex set 
obtained in \cite{T20} in the following theorem.

\begin{theorem}\cite{T20}\label{T:ThetaClassDaisyCube[T20]} 
Let $G$ be a daisy cube  whose vertex set is contained in $(\mathcal{B}^n, \le)$. 
Let $xy$ be an edge of $G$ whose two end vertices differ in exactly one position $i$ for some 
integer $i \in \{1, 2, \ldots, n\}$.
Then $F_{xy}=\{e \in E(G) \mid e \Theta xy\}$ is a peripheral $\Theta$-class, and
$G-F_{xy}$ has exactly two components $\langle W_{xy} \rangle$ and $\langle W_{yx} \rangle$ both of which are daisy cubes
such that one of them  contains all vertices of $G$ having $0$ at position $i$, 
and the other contains all vertices of $G$ having $1$ at position $i$.
Moreover, if $|W_{xy}| > |W_{yx}|$, 
then $\langle U_{xy} \rangle$ is a $\le$-subgraph of $G$, 
$W_{xy}$ is the subset of all vertices of $G$ having $0$ at position $i$, 
and $W_{yx}=U_{yx}$ is the subset of all vertices of $G$ having $1$ at position $i$.
\end{theorem}

\subsection{Decompositions of  resonance graphs} 
The concept of a reducible face decomposition of a plane bipartite graph $G$, 
denoted by $\mathrm{RFD}(G_1, G_2, \ldots, G_n)$, was introduced in \cite{ZZ00}.
Start from an edge $e$, add a path $P_1$ of odd length which has two end vertices in common with these of $e$
such that $e$ and $P_1$ form an even cycle $G_1$ surrounding a finite face $s_1$ of $G$,
proceed inductively to build a sequence of plane bipartite
graphs $G_i$ for $2 \le i \le n$ where $G_n=G$ as follows.
If $G_{i-1} =e + P_1+ \cdots + P_{i-1}$ has already been constructed, 
then $G_i=G_{i-1}+P_i$ can be obtained by adding the $i$th path $P_i$ of odd length
in the exterior of $G_{i-1}$ 
such that $P_i$ has exactly two end vertices in common with $G_{i-1}$,
$P_i$ and a part of the periphery of $G_{i-1}$ form an even cycle surrounding a finite face $s_i$ of $G$.
A sequence of finite faces $s_1, s_2, \ldots, s_n$ of $G$ associated with $\mathrm{RFD}(G_1, G_2, \ldots, G_n)$ is  generated 
during the above process.
It was shown \cite{ZZ00} that a plane bipartite graph $G$ with more than two vertices is elementary 
if and only if it has a reducible face decomposition.
If a plane elementary bipartite graph $G$ has a $\mathrm{RFD}(G_1, G_2, \ldots, G_n)$ associated with a sequence 
of finite faces $s_1, s_2, \ldots, s_n$, then it is easily seen that $s_i$ is a reducible face of $G_i$ for $2 \le i \le n$. 

It is well known  \cite{ZLS08} that the resonance graph of  a plane weakly elementary bipartite graph is a median graph.
Note that the resonance graph of an even cycle  is the one-edge graph.
Motivated from the structure properties of a median graph stated in Theorem \ref{T:Median},
it was proved in \cite{C18} that if $G$ is a plane elementary bipartite graph with more than two vertices,
then $R(G)$ can be constructed from the one-edge graph 
by a sequence of convex expansions with respect to a reducible face decomposition of $G$.
For a special case that $G$ has a $\mathrm{RFD}(G_1, G_2, \ldots, G_n)$ associated with a sequence 
of finite faces $s_1, s_2, \ldots, s_n $ such that  $s_i$ has common edges 
with exactly one other finite face $s_{\alpha(i)}$ in $G_i$ for $2 \le i \le n$,
then $R(G)$ can be constructed from the one-edge graph 
by a sequence of peripheral convex expansions with respect to $\mathrm{RFD}(G_1, G_2, \ldots, G_n)$.
It was further shown in \cite{C21} that the resonance graph of a plane elementary bipartite graph $G$ with more than two vertices
can be constructed from the one-edge graph 
by a sequence of peripheral convex expansions with respect to a reducible face decomposition of $G$
if and only if the infinite face of $G$ is forcing, that is, the removal of all vertices from the peripheral cycle of $G$
results in either an empty graph or a graph with an unique perfect matching.

\subsection{Subsets of perfect matchings of peripherally 2-colorable graphs} 

A cycle (or, a path) of a graph with a perfect matching $M$ is called \textit{$M$-alternating} 
if its edges are in and  out of $M$ alternately.  
A face (including the infinite face) of a plane bipartite graph with a perfect matching $M$  
is called \textit{$M$-resonant} if the periphery of the face is an $M$-alternating cycle.
A cycle that is the periphery of a face in a plane graph is called a \textit{facial cycle}.
A facial cycle of a finite face $s$ is denoted by $\partial s$. The facial cycle of the infinite face
is the periphery of $G$ and is denoted by $\partial G$.

Recall that a \textit{handle}  
of a graph is a path that has exactly two end vertices with degree larger than 2 in the graph \cite{CC13}. 
A  handle with more than one edge is called a \textit{nontrivial handle}. 
A handle of a plane graph   is called an \textit{exterior handle} 
if all of its edges are on the periphery of the plane graph, and called an \textit{interior handle}
if all of its edges are interior edges of the plane graph.
The concept of a handle was used in \cite{C18} as a key tool to partition the set of all perfect matchings of a plane bipartite graph
and obtain a decomposition structure of its resonance graph.

Let $G$ be a graph and $\mathcal{M}(G)$ be the set of all perfect matchings of $G$.
Let $P$ be an odd length handle of $G$.
By definition, if $P$ is nontrivial, then $P$ is $M$-alternating for any perfect matching $M$ of $G$.
It is easily seen that for any odd length handle $P$ and a perfect matching $M$ of $G$,
either $M$ contains all end edges of $P$ or $M$ does not contain any end edge of $P$. 
Two notations  $\mathcal{M}(G; P^{-})$ and $\mathcal{M}(G; P^{+})$ were  introduced in \cite{C18}.
Let $\mathcal{M}(G; P^{-})$ be the set of all perfect matchings $M$ of $G$ such that  
$M$ does not contain any end edge of  $P$.
Let $\mathcal{M}(G; P^{+})$ be the set of all perfect matchings $M$ of $G$ such that  
$M$ contains the end edges of  $P$.
Then $\mathcal{M}(G)$ is a disjoint union of $\mathcal{M}(G; P^{-})$ and $\mathcal{M}(G; P^{+})$.

 A plane elementary bipartite graph $G$ with more than two vertices  is called \textit{peripherally 2-colorable} 
if every vertex of $G$ has degree 2 or 3, 
vertices with degree 3 (if exist) are all exterior vertices of $G$, 
and $G$ can be properly 2-colored black and white so that any two adjacent vertices do not have the same color,
and vertices with degree 3 (if exist) are alternatively black and white along the clockwise orientation of the periphery of $G$. 

Let $G$ be a peripherally 2-colorable graph that is not an even cycle.
By definition, we can see that each exterior handle of $G$ has odd length, and any reducible face of $G$ 
has common edges with exactly one other finite face of $G$. 
Furthermore, by the induction proof based on a reducible face decomposition of $G$, it can be easily shown 
that each interior handle of $G$ has odd length, and any facial cycle of $G$ is a union of odd length handles that
are interior and exterior  alternately along the clockwise orientation of the facial cycle such that any two consecutive ones 
have exactly one vertex in common.

\vskip 0.1in
\begin{notations}\label{N}
Let $G$ be a peripherally 2-colorable graph that is not an even cycle.
Let $s$ be an arbitrary finite face of  $G$.
Let $\{P_{i}[s] \mid 1 \le i \le m\}$  (respectively, $\{J_{i}[s] \mid 1 \le i \le m\}$) 
be the set of exterior handles of $G$ (respectively, the set of interior handles of $G$)  
on the facial cycle $\partial s$ such that
$\partial s$ can be represented as 
$J_{1}[s] \cup P_{1}[s]  \cup \cdots \cup J_{m}[s] \cup P_{m}[s]$,
where  $J_{1}[s], P_{1}[s], \ldots, J_{m}[s], P_{m}[s]$ are listed in order along the clockwise orientation of $\partial s$,
and any two consecutive handles have exactly one  vertex in common.
\end{notations}
Notations \ref{N} introduced in the above
and the following Definition \ref{D:S-1-factors} will be applied throughout the next section without 
further illustrations.

\begin{definition}\label{D:S-1-factors}
Let $G$ be a  peripherally 2-colorable graph  that is not an even cycle.
Let $\mathcal{M}(G)$ be the set of all perfect matchings of $G$.
Let $i_0$ be an arbitrary integer in $\{1, \ldots,  m\}$. Then we can define the following subsets of $\mathcal{M}(G)$.
\begin{enumerate}
\item Let $\mathcal{M}(G; P^{-}_{i_0}[s])$ (respectively, $\mathcal{M}(G; J^{-}_{i_0}[s])$) be the set of perfect matchings $M$ of $G$ 
such that  $M$ does not contain any end edge of the exterior handle $P_{i_0}[s]$ of $G$ (respectively, the interior handle $J_{i_0}[s]$ of $G$). 

\item Let $\mathcal{M}(G; P^{+}_{i_0}[s])$  (respectively, $\mathcal{M}(G; J^{+}_{i_0}[s])$) be the set of perfect matchings $M$ of $G$
such that $M$ contains the end edges of the exterior handle $P_{i_0}[s]$ of $G$ (respectively, the interior handle $J_{i_0}[s]$ of $G$).

\item Let $\mathcal{M}(G; \mathcal{P}_G^{-}[s])$  (respectively, $\mathcal{M}(G; \mathcal{J}_G^{-}[s])$)
be the set of perfect matchings $M$ of $G$
such that $M$ does not contain any end edge of any exterior handle of $G$  (respectively, any interior handle of $G$)
on the facial cycle $\partial s$.

\item Let $\mathcal{M}(G; \mathcal{P}_G^{+}[s])$  (respectively, $\mathcal{M}(G; \mathcal{J}_G^{+}[s])$)
be the set of perfect matchings $M$ of $G$
such that $M$ contains end edges of all  exterior handles of $G$ (respectively, all interior handles of $G$)
on the facial cycle $\partial s$.

\item Let $\mathcal{M}(G; \mathcal{P}_G^{-}[s], \partial s)$ 
(respectively, $\mathcal{M}(G; \mathcal{J}_G^{-}[s], \partial s)$) 
be the set of perfect matchings $M$ of $G$ 
contained in $\mathcal{M}(G; \mathcal{P}_G^{-}[s])$ 
(respectively, $\mathcal{M}(G; \mathcal{J}_G^{-}[s])$)
such that $s$ is $M$-resonant.

\item Let $\mathcal{M}(G; \mathcal{P}_G^{+}[s], \partial s)$ 
(respectively, $\mathcal{M}(G; \mathcal{J}_G^{+}[s], \partial s)$)
be the set of perfect matchings $M$ of $G$ 
contained in $\mathcal{M}(G; \mathcal{P}_G^{+}[s])$
(respectively, $\mathcal{M}(G; \mathcal{J}_G^{+}[s])$)
such that $s$ is $M$-resonant.
\end{enumerate}
\end{definition}

\section{Main Results} 
\begin{lemma}\label{L:One-All-ExteriorHandles}
Let $G$ be a peripherally 2-colorable graph that is not an even cycle.
Let $s$ be an arbitrary finite face of  $G$. 
Let $i_0$ be an arbitrary integer in $\{1, \ldots,  m\}$.  Then

(i) $\mathcal{M}(G; P^{-}_{i_0}[s]) =  \mathcal{M}(G; \mathcal{P}_G^{-}[s])$. 

(ii)  $\mathcal{M}(G; P^{+}_{i_0}[s]) =  \mathcal{M}(G; \mathcal{P}_G^{+}[s])$.

(iii) $\mathcal{M}(G; \mathcal{P}_G^{+}[s])  
= \mathcal{M}(G; \mathcal{P}_G^{+}[s], \partial s)  
= \mathcal{M}(G; \mathcal{J}_G^{-}[s], \partial s)$.

(iv) $\mathcal{M}(G; \mathcal{J}_G^{+}[s]) 
= \mathcal{M}(G; \mathcal{J}_G^{+}[s], \partial s) 
= \mathcal{M}(G; \mathcal{P}_G^{-}[s], \partial s)$.

\end{lemma}
\proof  Recall that $\partial s$ is represented as 
$J_{1}[s] \cup P_{1}[s]  \cup \cdots \cup J_{m}[s] \cup P_{m}[s]$  
where handles $J_{1}[s], P_{1}[s], \ldots, J_{m}[s], P_{m}[s]$ are listed in order along the clockwise orientation of $\partial s$
such that any two consecutive ones have exactly one  vertex in common.

Let $i_0$ be an arbitrary integer in $\{1, \ldots,  m\}$.  
It is clear that $\mathcal{M}(G; \mathcal{P}_G^{-}[s]) \subseteq \mathcal{M}(G; P^{-}_{i_0}[s])$.
Suppose that $\mathcal{M}(G; P^{-}_{i_0}[s]) \neq  \mathcal{M}(G; \mathcal{P}_G^{-}[s])$. 
Then $G$ has a perfect matching $M$ such that $M \in \mathcal{M}(G; P^{-}_{i_0}[s])  \setminus \mathcal{M}(G; \mathcal{P}_G^{-}[s])$.
We observe that any handle of $G$ is $M$-alternating and with odd length.
Then $M$ either contains all end edges of a handle, or $M$ does not contain any end edge of a handle.
Without loss of generality, we can assume that $M \in \mathcal{M}(G; P^{-}_{i_0}[s]) \cap \mathcal{M}(G; P^{+}_{l}[s])$,
where $l=i_0+1$ if $i_0 \neq m$, and $l=1$ if $i_0 = m$. 
Let  $u_{i_0}$ and $v_{i_0}$ be two end vertices $P_{i_0}[s]$, 
and $u_{l}$ and $v_{l}$ be two end vertices $P_{l}[s]$ along the clockwise orientation of  $\partial s$.
Then the interior handle $J_{l}[s]$ of $G$ on $\partial s$ has two end vertices $v_{i_0}$ and $u_{l}$. 

Let $H$ be the subgraph obtained from $G$ by removing 
all internal vertices and edges from $P_{i}[s]$ for all $1 \le i \le m$.
Then $H$ has $m$ components $H_1, \ldots, H_{m}$,  
where $H_i$ contains the  interior handle  $J_{i}[s]$ of $G$ on $\partial s$ for $1 \le i \le m$.
Since $G$ is peripherally 2-colorable, it is easy to check that $H_i$ is peripherally 2-colorable for $1 \le i \le m$.

Let $H'$ be the subgraph obtained from $G$ by removing 
all internal vertices and edges from $P_{i_0}[s]$ and all vertices from $P_{l}[s]$.
Since $M \in \mathcal{M}(G; P^{-}_{i_0}[s]) \cap \mathcal{M}(G; P^{+}_{l}[s])$,  
the restriction of $M$ on $H'$ is a perfect matching of $H'$.
On the other hand, $H'$ has a component $H_{l} \setminus \{u_{l}\}$ with an odd number of vertices, and so 
$H'$ has no perfect matchings. This is a contradiction. 
Therefore, $\mathcal{M}(G; P^{-}_{i_0}[s]) =  \mathcal{M}(G; \mathcal{P}_G^{-}[s])$.
Similarly, we can show that  $\mathcal{M}(G; P^{+}_{i_0}[s]) =  \mathcal{M}(G; \mathcal{P}_G^{+}[s])$.

 By the structure of a peripherally 2-colorable graph $G$, 
if a perfect matching $M$ of $G$ contains end edges of $P_{i}[s]$ for all $1 \le i \le m$, then $\partial s$ is $M$-alternating,
that is, $s$ is $M$-resonant.
It follows that $M$  contains end edges of $P_{i}[s]$ for all $1 \le i \le m$  if and only if $s$ is $M$-resonant and 
 $M$  contains end edges of $P_{i}[s]$ for all $1 \le i \le m$. Furthermore, $s$ is $M$-resonant and 
 $M$  contains end edges of $P_{i}[s]$ for all $1 \le i \le m$ if and only if $s$ is $M$-resonant and 
$M$ contains no end edges of $J_{i}[s]$ for all $1 \le i \le m$.
Therefore, $\mathcal{M}(G; \mathcal{P}_G^{+}[s])  
= \mathcal{M}(G; \mathcal{P}_G^{+}[s], \partial s)  
= \mathcal{M}(G; \mathcal{J}_G^{-}[s], \partial s)$.

Similarly, $\mathcal{M}(G; \mathcal{J}_G^{+}[s]) 
= \mathcal{M}(G; \mathcal{J}_G^{+}[s], \partial s) 
= \mathcal{M}(G; \mathcal{P}_G^{-}[s], \partial s)$.
\qed\\

By Proposition 2.1 in \cite{T20},  each $\Theta$-class of a daisy cube is peripheral.
By Propositions 2.7-2.9 in \cite{T20}, any $\le$-subgraph of a daisy cube  is isomorphic to a daisy cube;
if $ab$ is an arbitrary edge of a daisy cube $G$ with a proper labelling such that vertex $a$ has $0$ (respectively, vertex $b$ has $1$) at position $i$, 
then $\langle U_{ab} \rangle$ is a $\le$-subgraph of $G$;
and  if $H$ is a $\le$-subgraph of a daisy cube $G$, then the $\le$-expansion of $G$ with respect to $H$ is a  daisy cube. 

By Lemma 3.4 and Corollary 3.6 in  \cite{BCTZ25},   if $G$ is a peripherally 2-colorable graph,
then $R(G)$ is a daisy cube and there is a 1--1 correspondence between the set of $\Theta$-classes of $R(G)$ and the set of finite faces of $G$
such that each finite face of  $G$ is the face-label of a unique $\Theta$-class of $R(G)$.
We are motivated to provide
a peripheral convex $\le$-decomposition structure of $R(G)$ that is a daisy cube with respect to an arbitrary finite face of $G$,
together with a proper labelling for the vertex set of $R(G)$.

\begin{theorem}\label{T:ArbitrayFiniteFace} 
Let $G$ be a peripherally 2-colorable graph that is not an even cycle.
Let $s$ be an arbitrary finite face of  $G$.
 Let $H$ be the subgraph of $G$ obtained by removing   
internal vertices and edges from all  exterior handles  of $G$ on $\partial s$.
Let $F$ be the set of edges in the resonance graph $R(G)$ with the face-label $s$.
Then

(i)  $F$ is a peripheral $\Theta$-class of $R(G)$ such that
the spanning subgraph $R(G)-F$ has exactly two components
$\langle \mathcal{M}(G; \mathcal{P}_G^{-}[s]) \rangle$ 
and $\langle  \mathcal{M}(G; \mathcal{P}_G^{+}[s], \partial s) \rangle$.

(ii)  $R(G)$ is a daisy cube whose vertex set has a proper labelling  in $\mathcal{B}^n$ 
such that the finite face $s=s_i$ of $G$ corresponds to the $i$-th position of the proper labelling, where $1 \le i \le n$.
Moreover,  $\langle \mathcal{M}(G; \mathcal{P}_G^{-}[s]) \rangle$  is a daisy cube 
that contains all vertices of $R(G)$ with a proper labelling having $0$ at position $i$, 
and $\langle  \mathcal{M}(G; \mathcal{P}_G^{+}[s], \partial s) \rangle$ is a daisy cube 
that contains all vertices of $R(G)$ with a proper labelling having $1$ at the same position $i$. 

(iii) The resonance graph $R(H)$ is a daisy cube and isomorphic to $ \langle \mathcal{M}(G; \mathcal{P}_G^{-}[s]) \rangle$.
A proper labelling for each vertex of $R(H)$ can be obtained by deleting $0$ at position $i$  from the proper labelling 
for the corresponding vertex of $\langle \mathcal{M}(G; \mathcal{P}_G^{-}[s]) \rangle$.

(iv) $\langle \mathcal{M}(G; \mathcal{P}_G^{-}[s], \partial s) \rangle
=\langle \mathcal{M}(G; \mathcal{J}_G^{+}[s]) \rangle \cong 
\langle \mathcal{M}(H; \mathcal{J}_G^{+}[s]) \rangle$,
and $F$ is a matching defining an isomorphism between two daisy cubes
$\langle \mathcal{M}(G; \mathcal{J}_G^{+}[s]) \rangle$ 
and $\langle \mathcal{M}(G; \mathcal{P}_G^{+}[s], \partial s) \rangle$.

(v) $R(G)$   can be obtained from $R(H)$ by a peripheral convex $\le$-expansion with respect to 
$\langle \mathcal{M}(H; \mathcal{J}_G^{+}[s]) \rangle$.
 \end{theorem} 

\proof $(i)$ If  $G$ is peripherally 2-colorable with $n$ finite faces,
then by Corollary 3.6 in \cite{BCTZ25},  
the resonance graph $R(G)$ is a daisy cube with $\mathrm{idim}(R(G)) = n$.

By Lemma 3.4 in \cite{BCTZ25},
there is a 1--1 correspondence between the set of finite faces of  $G$
and the set of $\Theta$-classes of $R(G)$ such that each finite face of $G$
is the face-label of all edges of a unique $\Theta$-class of $R(G)$. 
Let $F$ be the set of all edges in $R(G)$ with the face-label $s$.
It follows that $F$ is a $\Theta$-class of $R(G)$.
Assume that $xy$ is an edge in $F$.  Then $F=F_{xy}$ where $x \oplus y=\partial s$ since $xy$ has the face-label $s$.

Since $R(G)$ is a daisy cube,
by Theorem \ref{T:ThetaClassDaisyCube[T20]}, $F$ is peripheral,
and an edge is contained in $F$ if and only if the binary strings of its two end vertices from a proper labelling
differ in exactly one position $i$  for some integer $i \in \{1, 2, \ldots, n\}$.  
It follows that the vertex set of $R(G)$ has a proper labelling in $\mathcal{B}^n$
such that the finite face $s=s_i$ of $G$ corresponds to a unique position $i$ of the proper labelling, where $1 \le i \le n$.

Let $\mathcal{M}(G)$ be the set of all perfect matchings of $G$. 
Let $P_{i_0}[s]$ be an arbitrary exterior handle of $G$ on the facial cycle $s$.
Then for any $M \in \mathcal{M}(G)$, $P_{i_0}[s]$ is an $M$-alternating odd length handle
such that either $M$ contains all end edge(s) of $P_{i_0}[s]$ or $M$ contains no end edge(s) of $P_{i_0}[s]$.
Hence,  
$\mathcal{M}(G)=\mathcal{M}(G; P^{-}_{i_0}[s]) \cup \mathcal{M}(G; P^{+}_{i_0}[s])$.
By Lemma \ref{L:One-All-ExteriorHandles} $(i)$ - $(iii)$, we have $\mathcal{M}(G; P^{-}_{i_0}[s])=\mathcal{M}(G; \mathcal{P}_G^{-}[s])$
and $\mathcal{M}(G; P^{+}_{i_0}[s])=\mathcal{M}(G; \mathcal{P}_G^{+}[s])=\mathcal{M}(G; \mathcal{P}_G^{+}[s], \partial s)$.
Then $\mathcal{M}(G)= \mathcal{M}(G; \mathcal{P}_G^{-}[s]) \cup 
\mathcal{M}(G; \mathcal{P}_G^{+}[s], \partial s)$.
Recall that any edge of $R(G)$ with the face-label $s$ 
have the property that the symmetric difference of its two end vertices is the facial cycle $\partial s$.
Then any edge of $R(G)$ with the face-label $s$ must have one end vertex in $\langle \mathcal{M}(G; \mathcal{P}_G^{-}[s]) \rangle$
and the other end vertex in $\langle  \mathcal{M}(G; \mathcal{P}_G^{+}[s], \partial s) \rangle$. 
It follows that the set $F=F_{xy}$ of all edges in $R(G)$ with the face-label $s$ is the set of all edges between
$\langle \mathcal{M}(G; \mathcal{P}_G^{-}[s]) \rangle$ and
$\langle  \mathcal{M}(G; \mathcal{P}_G^{+}[s], \partial s) \rangle$. 
Hence, the spanning subgraph $R(G)-F$ has exactly two components
$\langle \mathcal{M}(G; \mathcal{P}_G^{-}[s]) \rangle$ 
and $\langle  \mathcal{M}(G; \mathcal{P}_G^{+}[s], \partial s) \rangle$.

$(ii)$ By definition, any peripherally 2-colorable graph $G$ is
plane elementary bipartite. By Corollary 3.3 in \cite{C18} (Equivalently, Theorem 3.8 in \cite{ZYY14}), 
$R(G)$ cannot be a nontrivial Cartesian product for any plane elementary bipartite graph $G$.
Note that $F$ is a matching between $\mathcal{M}(G; \mathcal{P}_G^{-}[s], \partial s)$ and
$\mathcal{M}(G; \mathcal{P}_G^{+}[s], \partial s)$.
Then $|\mathcal{M}(G; \mathcal{P}_G^{-}[s], \partial s)|= | \mathcal{M}(G; \mathcal{P}_G^{+}[s], \partial s)|$,
and the induced subgraph $\langle \mathcal{M}(G; \mathcal{P}_G^{-}[s], \partial s) \cup \mathcal{M}(G; \mathcal{P}_G^{+}[s], \partial s) \rangle$
of $R(G)$ is isomorphic to a Cartesian product $\langle \mathcal{M}(G; \mathcal{P}_G^{-}[s], \partial s) \rangle \Box K_2$.
By definitions, we know that the vertex set of $R(G)$ is 
$\mathcal{M}(G)=\mathcal{M}(G; \mathcal{P}_G^{-}[s]) \cup  \mathcal{M}(G; \mathcal{P}_G^{+}[s], \partial s)$,
and $\mathcal{M}(G; \mathcal{P}_G^{-}[s], \partial s)$
is contained in $\mathcal{M}(G; \mathcal{P}_G^{-}[s])$.
Since $R(G)$ cannot be a nontrivial Cartesian product, it follows 
that $\mathcal{M}(G; \mathcal{P}_G^{-}[s], \partial s)$
is properly contained in $\mathcal{M}(G; \mathcal{P}_G^{-}[s])$.
Hence, $|\mathcal{M}(G; \mathcal{P}_G^{-}[s])|  \neq  |  \mathcal{M}(G; \mathcal{P}_G^{+}[s], \partial s)|$,
and so $|\mathcal{M}(G; \mathcal{P}_G^{-}[s])| > |  \mathcal{M}(G; \mathcal{P}_G^{+}[s], \partial s)|$.  
Recall that $xy \in F$ with $x \oplus y=\partial s$. 
Without loss of generality, we can assume that $x \in \mathcal{M}(G; \mathcal{P}_G^{-}[s])$
and $y \in \mathcal{M}(G; \mathcal{P}_G^{+}[s], \partial s)$.
By Theorem \ref{T:ThetaClassDaisyCube[T20]}, 
both induced subgraphs $\langle W_{xy} \rangle=\langle \mathcal{M}(G; \mathcal{P}_G^{-}[s]) \rangle$ 
and $\langle W_{yx} \rangle=\langle  \mathcal{M}(G; \mathcal{P}_G^{+}[s], \partial s) \rangle$
are daisy cubes with a proper labelling  in $\mathcal{B}^n$ for their vertex sets 
such that $W_{xy}=\mathcal{M}(G; \mathcal{P}_G^{-}[s])$ contains all vertices of $R(G)$ having $0$ at position $i$ 
corresponding to the finite face $s(=s_i)$, 
and $W_{yx}=U_{yx}=\mathcal{M}(G; \mathcal{P}_G^{+}[s], \partial s)$  
contains all vertices of $R(G)$ having $1$ at  the same position $i$.

$(iii)$  Let $H$ be the subgraph of $G$ obtained by removing   
internal vertices and edges from all  exterior handles  of $G$ on $\partial s$.
Then any perfect matching of $H$ can be extended to a unique perfect matching of $G$ in $\mathcal{M}(G; \mathcal{P}_G^{-}[s])$.
On the other hand, the restriction of any perfect matching  of $G$ in $\mathcal{M}(G; \mathcal{P}_G^{-}[s])$ on $H$ 
is a perfect matching of $H$. Let $m_1$ and $m_2$ be two perfect matchings of $H$, and  $M_1$ and $M_2$ be the perfect matchings 
of $G$ in $\mathcal{M}(G; \mathcal{P}_G^{-}[s])$ extended from $m_1$ and $m_2$, respectively.

Then $m_1$ and $m_2$ are adjacent in $R(H)$ 

$\Longleftrightarrow$  $m_1 \oplus m_2 =\partial s_j$ where $s_j$ is a finite face of $G$ different from $s(=s_i)$

$\Longleftrightarrow$ $M_1 \oplus M_2=\partial s_j$ where $s_j$ is a finite face of $G$ different from $s(=s_i)$ 

$\Longleftrightarrow$  $M_1$ and $M_2$ are adjacent in $R(G)$.

Therefore, $R(H)$ is  isomorphic to   $\langle \mathcal{M}(G; \mathcal{P}_G^{-}[s]) \rangle$.
This implies that $R(H)$ is a daisy cube, and a proper labelling  in $\mathcal{B}^{n-1}$ for each vertex of $R(H)$ can be obtained 
by removing $0$ from the $i$-th position of the proper labelling  in $\mathcal{B}^n$
for the corresponding vertex of  $\langle \mathcal{M}(G; \mathcal{P}_G^{-}[s]) \rangle$.

$(iv)$ By Lemma \ref{L:One-All-ExteriorHandles},
$\mathcal{M}(G; \mathcal{P}_G^{-}[s], \partial s)=\mathcal{M}(G; \mathcal{J}_G^{+}[s])$.
Then $\langle \mathcal{M}(G; \mathcal{P}_G^{-}[s], \partial s) \rangle=
\langle \mathcal{M}(G; \mathcal{J}_G^{+}[s]) \rangle$.
Moreover, by the structure of a peripherally 2-colorable graph $G$,
the restriction of any perfect matching of $G$ in $\mathcal{M}(G; \mathcal{J}_G^{+}[s])$  on $H$ is a perfect matching of $H$
in $\mathcal{M}(H; \mathcal{J}_G^{+}[s])$. On the other hand, any perfect matching of $H$
in $\mathcal{M}(H; \mathcal{J}_G^{+}[s])$ can be extended to a unique perfect matching of $G$
in $\mathcal{M}(G; \mathcal{J}_G^{+}[s])$. 
In a way similiar to the proof for $(iii)$, we can show that 
$\langle \mathcal{M}(G; \mathcal{J}_G^{+}[s]) \rangle$ is isomorphic to $\langle \mathcal{M}(H; \mathcal{J}_G^{+}[s]) \rangle$.
Hence, $\langle \mathcal{M}(G; \mathcal{P}_G^{-}[s], \partial s) \rangle=
\langle \mathcal{M}(G; \mathcal{J}_G^{+}[s]) \rangle \cong \langle \mathcal{M}(H; \mathcal{J}_G^{+}[s]) \rangle$.

By Theorem 3.1 in \cite{ZLS08}, $R(G)$ is a median graph since $G$ is a plane elementary bipartite graph. 
By Theorem \ref{T:Median} (i), 
$F$ is a matching defining an isomorphism between 
$\langle U_{xy} \rangle = \langle \mathcal{M}(G; \mathcal{P}_G^{-}[s], \partial s) \rangle=
\langle \mathcal{M}(G; \mathcal{J}_G^{+}[s]) \rangle$ 
and $\langle U_{yx} \rangle =\langle \mathcal{M}(G; \mathcal{P}_G^{+}[s], \partial s) \rangle$.
We have shown that $\langle U_{yx} \rangle=\langle W_{yx} \rangle=\langle \mathcal{M}(G; \mathcal{P}_G^{+}[s], \partial s) \rangle$ 
is a daisy cube in $(ii)$.
By Theorem \ref{T:ThetaClassDaisyCube[T20]}, 
$\langle U_{xy} \rangle = \langle \mathcal{M}(G; \mathcal{J}_G^{+}[s]) \rangle$ is a $\le$-subgraph of the daisy cube $R(G)$.
By  Proposition 2.7 in \cite{T20}, any $\le$-subgraph of a daisy cube  is isomorphic to a daisy cube.
So, $ \langle \mathcal{M}(G; \mathcal{J}_G^{+}[s]) \rangle$ is a daisy cube. 

$(v)$ By Theorem \ref{T:Median} (ii),  $\langle U_{xy} \rangle=\langle \mathcal{M}(G; \mathcal{J}_G^{+}[s]) \rangle$ is
 convex in $\langle W_{xy} \rangle=\langle \mathcal{M}(G; \mathcal{P}_G^{-}[s]) \rangle$.
 We have shown that $\langle \mathcal{M}(G; \mathcal{J}_G^{+}[s]) \rangle 
\cong \langle \mathcal{M}(H; \mathcal{J}_G^{+}[s]) \rangle$
and $\langle \mathcal{M}(G; \mathcal{P}_G^{-}[s]) \rangle \cong R(H)$.
Then $\langle \mathcal{M}(H; \mathcal{J}_G^{+}[s]) \rangle$
is convex in $R(H)$.
Recall that $\langle \mathcal{M}(G; \mathcal{J}_G^{+}[s]) \rangle$ is a $\le$-subgraph of $R(G)$.
By the above proper labelling for the vertex sets of $R(H)$ and $R(G)$ respectively, we have that 
$\langle  \mathcal{M}(H; \mathcal{J}_G^{+}[s])\rangle$ is a $\le$-subgraph of $R(H)$. 
Therefore,  $\langle  \mathcal{M}(H; \mathcal{J}_G^{+}[s]) \rangle$ is a convex $\le$-subgraph  of $R(H)$.
By $(i)$ - $(iv)$ and the definitions of expansions, 
we can see that $R(G)$ can be obtained from $R(H)$ by a peripheral convex $\le$-expansion with respect to 
 $\langle  \mathcal{M}(H; \mathcal{J}_G^{+}[s]) \rangle$.
\qed\\

Note that any reducible face of a peripherally 2-colorable graph $G$ has common edges with exactly one other finite face of $G$,
and so has a unique exterior handle and a unique interior handle on its facial cycle.
Next, we provide a peripheral convex $\le$-decomposition structure of $R(G)$ with respect to a reducible face of $G$ and having more information
on the proper labelling for the vertex set of $R(G)$ as a daisy cube.

\begin{corollary}\label{C:ReducibleFace} 
Let $G$ be a peripherally 2-colorable graph that is not an even cycle.
Let  $s_n$ be a reducible face of $G$.
Then $s_n$ has common edges with exactly one other finite face $s_{\alpha(n)}$ of $G$.
Let $P[s_n]$ be the unique exterior handle of $G$ on $\partial s_n$,
and $J[s_n]$ be the unique interior handle  of $G$ on $\partial s_n$. 
Let $H$ be the subgraph of $G$ obtained by removing internal vertices and edges of $P[s_n]$.
Let $F_n$ be the set of edges in the resonance graph $R(G)$ with the face-label $s_n$. Then  

(i) $F_n$ is a  peripheral  $\Theta$-class of $R(G)$ and  the spanning subgraph $R(G)-F_n$ has exactly two components
$\langle \mathcal{M}(G; P^{-}[s_n]) \rangle$ and $\langle \mathcal{M}(G; P^{+}[s_n], \partial s_n) \rangle$.

(ii) $R(G)$ is a daisy cube whose vertex set has a proper labelling  in $\mathcal{B}^n$ 
such that $s_n$ corresponds to the $n$-th position of the proper labelling.
Moreover, $\langle \mathcal{M}(G; P^{-}[s_n]) \rangle$ is a daisy cube that contains all vertices of $R(G)$ with a proper labelling having $0$ at position $n$, 
and $ \langle \mathcal{M}(G; P^{+}[s_n], \partial s_n)\rangle$ is a daisy cube that contains all vertices of $R(G)$ with a proper labelling having $1$ at  position $n$.

(iii) The resonance graph $R(H)$ is a daisy cube and isomorphic to  $\langle \mathcal{M}(G; P^{-}[s_n]) \rangle$.  A proper labelling for each vertex
of $R(H)$ can be obtained by deleting $0$ at position $n$  from the proper labelling 
for the corresponding vertex of $\langle \mathcal{M}(G; P^{-}[s_n])  \rangle$.

(iv) $\langle \mathcal{M}(G; P^{-}[s_n], \partial s_n) \rangle=\langle \mathcal{M}(G; J^{+}[s_n]) \rangle \cong \langle \mathcal{M}(H; J^{+}[s_n]) \rangle$, 
and $F_n$ is a matching defining an isomorphism between two daisy cubes
$\langle \mathcal{M}(G; J^{+}[s_n]) \rangle$ and $\langle  \mathcal{M}(G; P^{+}[s_n], \partial s_n) \rangle$.

(v) $R(G)$ can be obtained from $R(H)$ by a peripheral convex $\le$-expansion with respect to
$\langle \mathcal{M}(H; J^{+}[s_n]) \rangle$.

(vi) 
$\langle \mathcal{M}(H; J^{+}[s_n]) \rangle$ is a daisy cube that  contains all vertices of $R(H)$ with a proper labelling 
having $0$ at position $\alpha(n)$ corresponding to the finite face $s_{\alpha(n)}$,  and
$\langle  \mathcal{M}(G; J^{+}[s_n]) \rangle$  is a daisy cube that  contains all vertices of $R(G)$ with a proper labelling 
having $0$ at both position $\alpha(n)$ and position $n$.
 
\end{corollary} 

\proof  If $s_n$ is a reducible face of $G$, 
then the facial cycle $\partial s_n$ can be represented as $P[s_n] \cup J[s_n]$ 
where $P[s_n]$ is the unique exterior handle of $G$  on $\partial s_n$,
and $J[s_n]$ is the unique interior handle of $G$ on $\partial s_n$
such that $P[s_n]$ and $J[s_n]$ have two end vertices in common.
By Theorem \ref{T:ArbitrayFiniteFace}, the statements from  $(i)$ up to $(v)$  hold true immediately.
It remains to show that $(vi)$ holds true.
By the structure of a peripherally 2-colorable graph $G$ and the definition of a reducible face $s_n$ of $G$,
we can see that $s_n$ has common edges with exactly one other finite face $s_{\alpha(n)}$ of $G$.
Moreover, $\partial s_n \cap \partial s_{\alpha(n)}=J[s_n]$.

Let $H$ be the subgraph of $G$ obtained by removing internal vertices and edges of $P[s_n]$.
Then $H$ is a peripherally 2-colorable graph and $s_{\alpha(n)}$ is also a finite face of  $H$. 
Let $\{P_{H, i}[s_{\alpha(n)}] \mid 1 \le i \le m'\}$ be the set of exterior handles of $H$ 
on the facial cycle $\partial s_{\alpha(n)}$,
and  $\{J_{H, i}[s_{\alpha(n)}] \mid 1 \le i \le m'\}$ be the set of interior handles of $H$ 
on the facial cycle $\partial s_{\alpha(n)}$. 
We can assume that $\partial s_{\alpha(n)}$ is represented as 
$J_{H, 1}[s] \cup P_{H, 1}[s]  \cup  \ldots \cup J_{H, m'}[s] \cup P_{H, m'}[s]$ where 
$J_{H, 1}[s], P_{H, 1}[s], \ldots, J_{H, m'}[s], P_{H, m'}[s]$ are listed in order along the clockwise orientation of $\partial s_{\alpha(n)}$
such that any two consecutive ones have exactly one vertex in common.

We observe that  $J[s_n]$ is an interior handle of $G$ on the facial cycle $\partial s_{\alpha(n)}$,
and  contained in one  exterior handle of $H$ on the facial cycle  $\partial s_{\alpha(n)}$.
Without loss of generality, we can assume that  $J[s_n]$ is  contained in $P_{H, m'}[s_{\alpha(n)}]$.
Recall that the vertex degree of $G$ is at most 3. 
Then by the definition of a handle, any end vertex of a handle of a peripherally 2-colorable graph $H$ or $G$
has vertex degree 3 in both $H$ and $G$. It follows that 
any end vertex of $J[s_n]$ cannot be in common with any end vertex of $P_{H, m'}[s_{\alpha(n)}]$.

Now for $1 \le i \le m'-1$,  exterior handles $P_{H, i}[s_{\alpha(n)}]$ of $H$ on $\partial s_{\alpha(n)}$ are also
exterior handles of $G$ on $\partial s_{\alpha(n)}$; and
for $1 \le i \le m'$,  interior handles $J_{H, i}[s_{\alpha(n)}]$ of $H$ on $\partial s_{\alpha(n)}$ are also
interior handles of $G$ on $\partial s_{\alpha(n)}$. 

When $i=m'$, $P_{H, m'}[s_{\alpha(n)}]$
is a union of one interior handle $J[s_n]$ of $G$ and two exterior handles of $G$ on $\partial s_{\alpha(n)}$.
Let $P_{G,m'}[s_{\alpha(n)}]$ and $P_{G,m'+1}[s_{\alpha(n)}]$ be two such  exterior handles of $G$ on $\partial s_{\alpha(n)}$.
Then we can write $P_{H, m'}[s_{\alpha(n)}]=P_{G,m'}[s_{\alpha(n)}] \cup J[s_n] \cup P_{G,m'+1}[s_{\alpha(n)}]$.
Note that $P_{G,m'}[s_{\alpha(n)}]$, $J[s_n]$, and $P_{G,m'+1}[s_{\alpha(n)}]$
are odd length handles of $G$ such that any two consecutive ones have exactly one vertex in common. 
It follows that a perfect matching of $H$ is contained in $\mathcal{M}(H; J^{+}[s_n])$ 
if and only if it is contained in $\mathcal{M}(H; P^{-}_{H,m'}[s_{\alpha(n)}])$.
Hence, $\mathcal{M}(H; J^{+}[s_n])=\mathcal{M}(H; P^{-}_{H,m'}[s_{\alpha(n)}])$.

Let $\mathcal{M}(H; \mathcal{P}^{-}_{H}[s_{\alpha(n)}])$ be the set of perfect matchings $M$ of $H$
such that $M$ does not contain any end edge of any exterior handle of $H$ on the facial cycle $s_{\alpha(n)}$. 
By Lemma \ref{L:One-All-ExteriorHandles} $(i)$, 
$\mathcal{M}(H; P^{-}_{H,m'}[s_{\alpha(n)}])=\mathcal{M}(H; \mathcal{P}^{-}_{H}[s_{\alpha(n)}])$.
So, $\mathcal{M}(H; J^{+}[s_n])=\mathcal{M}(H; \mathcal{P}^{-}_{H}[s_{\alpha(n)}])$.

By Theorem \ref{T:ArbitrayFiniteFace} $(ii)$,  $\langle \mathcal{M}(H; \mathcal{P}^{-}_{H}[s_{\alpha(n)}])\rangle$  is a daisy cube 
that contains all vertices of $R(H)$ with a proper labelling having $0$ at position $\alpha(n)$  corresponding to the finite face $s_{\alpha(n)}$.
So, $\langle \mathcal{M}(H; J^{+}[s_n]) \rangle$  is a daisy cube 
that contains all vertices of $R(H)$ with a proper labelling having $0$ at position $\alpha(n)$.

By $(iv)$, $\langle \mathcal{M}(H; J^{+}[s_n]) \rangle \cong \langle \mathcal{M}(G; J^{+}[s_n]) \rangle
=\langle \mathcal{M}(G; P^{-}[s_n], \partial s_n) \rangle 
\subseteq \langle \mathcal{M}(G; P^{-}[s_n])  \rangle$.  By $(ii)$,
$\langle \mathcal{M}(G; P^{-}[s_n])  \rangle$ is a daisy cube that contains all vertices of $R(G)$ 
with a proper labelling having $0$ at position $n$.
By $(iii)$, $R(H)$ is isomorphic to $\langle \mathcal{M}(G; P^{-}[s_n])  \rangle$,
and a proper labelling for each vertex
of $R(H)$ can be obtained by deleting $0$ at position $n$  from the proper labelling 
of the corresponding vertex of $\langle \mathcal{M}(G; P^{-}[s_n])  \rangle$.
Hence, $\langle  \mathcal{M}(G; J^{+}[s_n]) \rangle$ is a daisy cube that contains all vertices of $R(G)$ 
with a proper labelling having  $0$ at both position $\alpha(n)$ and position $n$.
\qed\\

By \cite{BCTZ25},
the resonance graph $R(G)$ of a plane elementary bipartite graph $G$ with more than two vertices
is a daisy cube if and only if $G$ is peripherally 2-colorable.
Assume that $G$ is a  peripherally 2-colorable graph with a 
$\mathrm{RFD}(G_1, G_2, \ldots, G_n)$ associated 
with a sequence of finite faces $s_1, s_2, \ldots, s_n$ where each $s_i$ is a reducible face of $G_i$ 
for $2 \le i \le n$.
Then each $G_i$ is  peripherally 2-colorable for $1 \le i \le n$, 
and each $s_i$ has common edges with exactly one other finite face $s_{\alpha(i)}$ in $G_i$ for $2 \le i \le n$.
Applying Corollary \ref{C:ReducibleFace} repeatedly 
for the sequence of  peripherally 2-colorable graphs $G_i$ 
with a reducible face $s_i$ for $2 \le i \le n$,   
we  obtain the following corollary and an algorithm to generate a proper labelling for the vertex set of
the resonance graph $R(G)$ as a daisy cube with respect to a reducible face decomposition of $G$.

\begin{corollary}\label{C:pces-Daisycube}
Let $G$ be a plane elementary bipartite with more than two vertices. 
Then $R(G)$ is a daisy cube if and only if 
it can be obtained from the one-edge graph by a sequence of
peripheral convex $\le$-expansions with respect to a reducible face decomposition of $G$.
\end{corollary}

\begin{algorithm}\label{A:DaisyCube}
{\bf Input:} 
A peripherally 2-colorable graph $G$ with a $\mathrm{RFD}(G_1,G_2, \ldots ,G_n)$
 associated with a sequence of finite faces $s_1, s_2, \ldots, s_n$ where $n$ is a positive integer,
 and $s_i$ is a reducible face of $G_i$ for $2 \le i \le n$.

{\bf Output:} A  binary coding $\mathcal{L}_n$  for the vertex set of the resonance graph $R(G)$ as a daisy cube. 

Step 0. $i := 1$, $L_i := \{0,1\}$.

Step 1. If $i = n$, stop.

Step 2. Assume that $s_{i+1}$ has common edges with $s_{\alpha(i+1)}$, where $\alpha(i+1)$ is an integer in
$\{1, \ldots, i\}$, and $(x)_{\alpha(i+1)}$ is the $\alpha(i+1)$-th digit of a binary string $x$ of length $i$.
Set $L_{i+1} := \{x0 : x \in L_i\}  \cup  \{x1 : x \in L_i \  and \  (x)_{\alpha(i+1)} =0\}$.

Step 3. $i:=i+1$, go to step 1.
\end{algorithm}

See Figure  \ref{A-DaisyCube} for the illustration of Algorithm \ref{A:DaisyCube}.

\begin{figure}[h!]
\begin{center}
\begin{tikzpicture}[scale=0.5]
\tikzset{vertex/.style={draw=black, circle, minimum size=4pt, scale=0.5}}
 \draw [color=black, mark=*] plot[samples at={-150,-90,...,150,210},  variable=\x] 
  (\x:1);
    \node [vertex] at (-150:1){};
  \node  [vertex, fill] at (-90:1) {};
  \node [vertex] at (-30:1) {} ;
  \node [vertex, fill] at (30:1) {} ;
  \node  [vertex] at (90:1) {};
  \node [vertex, fill] at (150:1) {};
  \node [blue]{$s_1$};

 \begin{scope}[xshift=0.866 cm,  yshift=-1.5cm]
    \draw [color=black,mark=*] plot[samples at={-150,-90,...,210},variable=\x] 
  (\x:1);
      \node [vertex] at (-150:1){};
  \node  [vertex, fill] at (-90:1) {};
  \node [vertex] at (-30:1) {} ;
  \node [vertex, fill] at (30:1) {} ;
  \node  [vertex] at (90:1) {};
  \node [vertex, fill] at (150:1) {};
  \node[blue]{$s_2$};
  \end{scope}
    \begin{scope}[xshift=2.598 cm,  yshift=-1.5cm]
    \draw [color=black,mark=*] plot[samples at={-150,-90,...,210},variable=\x] 
  (\x:1);
      \node [vertex] at (-150:1){};
  \node  [vertex, fill] at (-90:1) {};
  \node [vertex] at (-30:1) {} ;
  \node [vertex, fill] at (30:1) {} ;
  \node  [vertex] at (90:1) {};
  \node [vertex, fill] at (150:1) {};
   \node [blue]{$s_3$};
  \end{scope}
      \begin{scope}[xshift=3.464 cm]
    \draw [color=black, mark=*] plot[samples at={-150,-90,...,210},variable=\x] 
  (\x:1);
      \node [vertex] at (-150:1){};
  \node  [vertex, fill] at (-90:1) {};
  \node [vertex] at (-30:1) {} ;
  \node [vertex, fill] at (30:1) {} ;
  \node  [vertex] at (90:1) {};
  \node [vertex, fill] at (150:1) {};
   \node [blue]{$s_4$};
  \end{scope}
  \begin{scope}[yshift=-3cm]
    \draw [color=black,mark=*] plot[samples at={-150,-90,...,210},variable=\x] 
  (\x:1);
      \node [vertex] at (-150:1){};
  \node  [vertex, fill] at (-90:1) {};
  \node [vertex] at (-30:1) {} ;
  \node [vertex, fill] at (30:1) {} ;
  \node  [vertex] at (90:1) {};
  \node [vertex, fill] at (150:1) {};
   \node [blue](s5){$s_5$};
   \node[draw=none, below of = s5,  xshift=0.7cm]{$G$};
   \end{scope}
 
\tikzstyle{every node}=[draw, circle, inner sep=2pt, minimum size=2pt, label distance=0.01pt, node distance=1cm]
\tikzset{el/.style = {inner sep=2pt, align=left}}
\tikzset{every label/.append style = {font=\tiny}}
\begin{scope}[xshift=9cm]
\node [ fill, label=right:$1$] (1) {};
\node [ fill, label=right:$\mathbf{0}$] (2) [below of=1] {};
\path[every node/.style={font=\sffamily\small}, inner sep=0.2pt, line width=1.pt]
(1) edge [line width= 3pt, color=blue] node[el, right, pos=0.5] {$s_1$} (2);
\node[draw=none, below of =1, yshift=-1cm]{$R(G_1)$};
\end{scope}

\begin{scope}[xshift=16cm]
\node [fill, label=right:$10$] (1) {};
\node [fill, label=right:$\mathbf{00}$] (2) [below of=1] {};
\node [fill, label=right:$01$] (3) [above of=2, xshift=-2cm] {};
\path[every node/.style={font=\sffamily\small}, inner sep=0.5pt, line width=1.pt]
(1) edge [color=blue] node[el, right, pos=0.5] {$s_1$} (2)
(2) edge [line width= 3pt, color=blue] node[el, right, pos=0.6] {$s_2$} (3);
\node[draw=none, below of =1, yshift=-1cm]{$R(G_2)$};
\end{scope}

\begin{scope}[xshift=23.5cm]
\node [fill, label=right:$100$] (1) {};
\node [fill, label=right:$\mathbf{000}$] (2) [below of=1] {};
\node [fill, label=right:$010$] (3) [above of=2, xshift=-2cm] {};
\node [fill, label=right:$001$] (4) [below of=1, xshift=1cm, yshift=1cm] {};
\node [fill, label=right:$101$] (5) [above of=2, xshift=1cm, yshift=1cm] {};
\path[every node/.style={font=\sffamily\small}, inner sep=0.5pt, line width=1.pt]
(1) edge [color=blue] node[el, midway, pos=0.5] {$s_1$} (2)
(4) edge [color=blue] node[el, right, pos=0.5] {$s_1$} (5)
(2) edge [color=blue] node[el, left, pos=0.5] {$s_2$} (3)
(1) edge [line width= 3pt, color=blue ] node[el, above, pos=0.4] {$s_3$} (5)
(2) edge [line width= 3pt, color=blue] node[el, right, pos=0.4] {$s_3$} (4);
\node[draw=none, below of =1, yshift=-1cm ]{$R(G_3)$};
\end{scope}
\end{tikzpicture}

\begin{tikzpicture}
\tikzstyle{every node}=[draw, circle, inner sep=2pt, minimum size=2pt, label distance=0.01pt, node distance=1cm]
\tikzset{el/.style = {inner sep=2pt, align=left}}
\tikzset{every label/.append style = {font=\tiny}}

\begin{scope}[xshift=2cm, yshift=-2.8cm]
\node [fill, label=right:$1000$] (1) {};
\node [fill, label=right:$\mathbf{0000}$] (2) [below of=1] {};
\node [fill, label=left:$0100$] (3) [above of=2, xshift=-2cm] {};
\node [fill, label=right:$0010$] (4)  [below of=1, xshift=1cm, yshift=1cm] {};
\node [fill, label=right:$1010$] (5) [above of=2,  xshift=1cm, yshift=1cm] {};
\node [fill, label=left:$0101$] (6) [above of=1,  xshift=-3cm] {};
\node [fill, label=left:$0001$] (7) [below of=1, xshift=-1cm, yshift=1cm] {};
\node [fill, label=right:$1001$] (8) [above of=7] {};
\path[every node/.style={font=\sffamily\small}, inner sep=0.5pt, line width=1.pt]
(1) edge [color=blue] node[el, midway, pos=0.5] {$s_1$} (2)
(4) edge [color=blue] node[el, right, pos=0.4] {$s_1$} (5)
(2) edge [color=blue] node[el, below, pos=0.5] {$s_2$} (3)
(1) edge [color=blue] node[el, left, pos=0.7] {$s_3$} (5)
(2) edge [color=blue] node[el, right, pos=0.5] {$s_3$} (4)
(6) edge [color=blue] node[el, right, pos=0.5] {$s_2$} (7)
(7) edge [color=blue] node[el, right, pos=0.3] {$s_1$} (8)
(1) edge [line width= 3pt, color=blue] node[el, right, pos=0.5] {$s_4$} (8)
(2) edge [line width= 3pt, color=blue] node[el, right, pos=0.7] {$s_4$} (7)
(3) edge [line width= 3pt, color=blue] node[el, left, pos=0.4] {$s_4$} (6);
\node[draw=none, below of =1, yshift=-1cm]{$R(G_4)$};
\end{scope}

\begin{scope}[xshift=9cm, yshift=-2.8cm]
\node [fill, label=left: $10000$] (1a) {};
\node [fill, label=right:$\mathbf{00000}$] (2a) [below of=1a] {};
\node [fill, label=left:$01000$] (3a) [above of=2a, xshift=-3cm] {};
\node [fill, label=right:$00100$] (4a)  [below of=1a, xshift=1.8cm, yshift=1cm] {};
\node [fill, label=left:$10100$] (5a) [above of=2a,  xshift=1.8cm, yshift=1cm] {};
\node [fill, label=right:$01010$] (6a) [above of=1a,  xshift=-4.5cm] {};
\node [fill, label=left:$00010$] (7a) [below of=1a, xshift=-1.5cm, yshift=1cm] {};
\node [fill, label=left:$10010$] (8a) [above of=7a] {};
\path[every node/.style={font=\sffamily\small}, inner sep=0.5pt, line width=1.pt]
(1a) edge [color=blue] node[el, midway, pos=0.4] {$s_1$} (2a)
(4a) edge [color=blue] node[el, midway, pos=0.7] {$s_1$} (5a)
(2a) edge [color=blue] node[el, below, pos=0.5] {$s_2$} (3a)
(1a) edge [color=blue] node[el, above, pos=0.5] {$s_3$} (5a)
(2a) edge [color=blue] node[el, right, pos=0.5] {$s_3$} (4a)
(1a) edge [color=blue] node[el, midway, pos=0.3] {$s_4$} (8a)
(2a) edge [color=blue] node[el, right, pos=0.7] {$s_4$} (7a)
(3a) edge [color=blue] node[el, left, pos=0.4] {$s_4$} (6a)
(6a) edge [color=blue] node[el, above, pos=0.5] {$s_2$} (7a)
(7a) edge [color=blue] node[el, left, pos=0.5] {$s_1$} (8a);
\end{scope}

\begin{scope}[xshift=9.5cm, yshift=-1.8cm]
\node [fill, label=above:$10001$, label distance=1pt] (1b) {};
\node [fill, label=right:$00001$] (2b) [below of=1b] {};
\node [fill, label=right:$00101$] (4b)  [below of=1b, xshift=1.8cm, yshift=1cm] {};
\node [fill, label=right:$10101$] (5b) [above of=2b,  xshift=1.8cm, yshift=1cm] {};
\node [fill, label=right:$00011$] (7b) [below of=1b, xshift=-1.5cm, yshift=1cm] {};
\node [fill, label=left:$10011$] (8b) [above of=7b] {};
\path[every node/.style={font=\sffamily\small}, inner sep=0.5pt, line width=1.pt]
(1b) edge [color=blue] node[el, midway, pos=0.4] {$s_1$} (2b)
(4b) edge [color=blue] node[el, right, pos=0.5] {$s_1$} (5b)
(1b) edge [color=blue] node[el, above, pos=0.5] {$s_3$} (5b)
(2b) edge [color=blue] node[el, below, pos=0.5] {$s_3$} (4b)
(1b) edge [color=blue] node[el, right, pos=0.7] {$s_4$} (8b)
(2b) edge [color=blue] node[el, right, pos=0.7] {$s_4$} (7b)
(7b) edge [color=blue] node[el, midway, pos=0.4] {$s_1$} (8b)
(1b) edge [line width= 3pt, color=blue] node[el, left, pos=0.2] { $s_5$} (1a)
(2b) edge [line width= 3pt, color=blue] node[el, right, pos=0.4] { $s_5$} (2a)
(4b) edge [line width= 3pt, color=blue] node[el, right, pos=0.6] { $s_5$} (4a)
(5b) edge [line width= 3pt, color=blue] node[el, left, pos=0.7] { $s_5$} (5a)
(7b) edge [line width= 3pt, color=blue] node[el, right, pos=0.7] { $s_5$} (7a)
(8b) edge [line width= 3pt, color=blue] node[el, left, pos=0.6] { $s_5$} (8a);
\node[draw=none, below of =1a,  yshift=-1cm]{$R(G)$};
\end{scope}
\end{tikzpicture}
\end{center}
\caption{\label{A-DaisyCube} An example for implementing Algorithm \ref{A:DaisyCube}.}
\end{figure}

By Lemma 3.4 in \cite{BCTZ25}, if $G$ is a peripherally 2-colorable graph with $n$ finite faces,
then  $\mathrm{idim}(R(G))=n$.  
\begin{corollary}
 Algorithm \ref{A:DaisyCube} can be applied to generate a proper labelling for all perfect matchings of a plane weakly elementary bipartite
 graph $G'$ whose each elementary component with more than two vertices is peripherally 2-colorable,
inducing an isometric embedding of $R(G')$ into $Q_{d}$ as a daisy cube, where $d=\mathrm{idim}(R(G'))$
is the sum of number of finites faces of all peripherally 2-colorable components of $G'$.
\end{corollary}
\proof By Theorem 4.3 in \cite{BCTZ25}, if $G'$ is a plane bipartite graph with a perfect matching, then
its resonance graph $R(G')$ is a daisy cube if and only if $G'$ is weakly elementary
such that each elementary component $G_i$  with more than two vertices is  peripherally 2-colorable where $1 \le i \le t$.
Moreover, $R(G')=R(G_1) \Box \cdots \Box R(G_t)$ where each $R(G_i)$ is a daisy cube for $1 \le i \le t$.
Assume that the number of finite faces of $G_i$ is  $n_i$ for all $1 \le i \le t$.  
By Lemma 3.4 in \cite{BCTZ25},  $\mathrm{idim}(R(G_i))=n_i$.
We can apply  Algorithm \ref{A:DaisyCube} to generate a proper labelling $\varphi_i: \mathcal{M}(G_i) \longrightarrow \mathcal{B}^{n_i}$ 
for all perfect matchings of $G_i$ for each $1 \le i \le t$. 
Let $\varphi=(\varphi_1\ldots \varphi_t)$ and $d=n_1+n_2 + \cdots +n_t$. Then 
$\varphi: \mathcal{M}(G') \longrightarrow \mathcal{B}^{d}$ is 
a proper labelling for all perfect matchings of $G'$.
By \cite{O08}, the isometric dimension of a Cartesian product graph is the summation 
of the isometric dimensions of all factors. 
It follows that $\mathrm{idim}(R(G'))= \sum\limits_{i=1}^{t}\mathrm{idim}(R(G_i))=\sum\limits_{i=1}^{t} n_i=d$. 
\qed\\

{\bf Remark} Algorithm  \ref{A:DaisyCube} is more efficient than Algorithm 2  in \cite{BCTZ24+} 
since we do not need to generate all maximal independent 
sets of the inner dual of a peripherally 2-colorable graph.

Algorithm  \ref{A:DaisyCube} is also presented more concisely comparing to Algorithm 1  in \cite{BCTZ24+},
although they have the same time complexity. This can be seen as follows.
Algorithm 1  in \cite{BCTZ24+} is exactly the first part of Algorithm 1 in \cite{BTZ23}
for the case of a 2-connected outerplane bipartite graph $G$  being ``regular'', while
Algorithm 1 in \cite{BTZ23} provides a binary coding for perfect matchings
of  a 2-connected outerplane bipartite graph $G$ based on whether $G$
is ``regular'' or not, and the latter is a generalization of the Algorithm in \cite{KVZG01} for catacondensed benzenoid systems.
By Corollary \ref{C:ReducibleFace}, we can see that there is no need to consider any adjacent triple of finite faces  
during the process of generating a binary coding for all perfect matchings of a peripherally 2-colorable graph. 
Hence, the line $i = \min \{ \ l \mid s_l \mbox{ is adjacent to } s_j \}$ of Algorithm 1 in \cite{BCTZ24+} is not needed at all. 
Moreover, it is not obvious that Algorithm 1  in \cite{BCTZ24+}  generates a binary coding with the described properties
for all perfect matchings of a peripherally 2-colorable graph.


In essence, we not only provide an efficient Algorithm  \ref{A:DaisyCube} 
to generate a proper labelling for all perfect matchings
of a peripherally 2-colorable graph $G$ which induces a daisy cube structure on $R(G)$,
but also prove Algorithm  \ref{A:DaisyCube} with Corollary \ref{C:ReducibleFace}.

\section{Two distinct binary codings}
It has been known that resonance graphs of various types of plane bipartite graphs display rich structures 
such as finite distributive lattices and daisy cubes.
Note that a finite distributive lattice and a daisy cube are distinct types of posets.  
A finite distributive lattice  has  a minimum element and a maximum element, while
a  daisy cube has a  minimum element  and a set of maximal elements.

Let $M$ be a perfect matching of a plane bipartite graph whose vertices are properly colored
white and black such that any two adjacent vertices receive different colors. 
Then an $M$-alternating cycle of the graph is called \textit{$M$-proper} 
(respectively, \textit{$M$-improper}) if
every edge of the cycle contained in $M$ goes from white to black vertices (respectively, from black
to white vertices) along the clockwise orientation of the cycle. 
An $M$-alternating path on a cycle in the graph is called \textit{$M$-proper} 
(respectively, \textit{$M$-improper}) if
every edge of the path contained in $M$ goes from white to black vertices (respectively, from black
to white vertices) along the clockwise orientation of the cycle. 

By \cite{LZ03}, the resonance graph $R(G)$ of any plane elementary bipartite graph $G$ has a finite distributive lattice structure 
with the minimum element being the unique perfect matching $M_{\hat{0}}$ of $G$ 
such that $G$ has no proper  $M_{\hat{0}}$-alternating cycles,
and  the maximum element being the unique  perfect matching $M_{\hat{1}}$ of $G$ 
such that $G$ has no improper  $M_{\hat{1}}$-alternating cycles.
Let $G$ be a plane elementary bipartite graph with $n$ finite faces and whose infinite face is forcing. 
By \cite{C24+, ZLS08}, there is a binary coding
$\mathcal{L}_n^{FDL}$ of length $n$ for all perfect matchings of $G$
which induces an isometric embedding of $R(G)$ into a hypercube $Q_n$ as a finite distributive lattice.

On the other hand, by \cite{BCTZ25}, 
resonance graphs  can be daisy cubes for exactly one special type of plane elementary bipartite graphs, 
namely peripherally 2-colorable graphs. Let $G$ be a peripherally 2-colorable graph with $n$ finite faces. 
By \cite{BCTZ25},  resonance graph $R(G)$ has a daisy cube structure 
with the minimum element $M_0$ being the unique perfect matching of $G$ such that every finite face of $G$ is $M_0$-resonant.
Moreover,  there is a binary coding $\mathcal{L}_n^{DC}$ for all perfect matchings of $G$ 
which induces an isometric embedding of $R(G)$ into a hypercube $Q_n$ as a daisy cube.

If $G$ is a peripherally 2-colorable graph with one finite face, then $G$ is an even cycle with exactly two perfect matchings, and $R(G)$ 
is the one-edge graph. Any binary coding of length 1 on $\mathcal{M}(G)$ assigns one perfect matching of $G$ as 0 and the other as 1.
 If $G$ is a peripherally 2-colorable graph with at least two finite faces,
then the infinite face of $G$ is forcing, and Algorithm 3.3 provided in \cite{C24+} generates a binary coding for all 
perfect matchings of $G$ inducing a distinct isometric embedding of $R(G)$  into $Q_n$ as a finite distributive lattice,
where $n=\mathrm{idim}(R(G))$ is the number of finite faces of $G$. 
Note that Algorithm 3.3   in \cite{C24+}  is a generalization of Algorithm 5.1 in \cite{ZLS08} from catacondensed
benzenoid graphs and 2-connected outerplane bipartite graphs to all plane elementary bipartite graphs 
whose infinite face is forcing.
Moreover,  if $G'$ is a plane weakly elementary bipartite graph
whose each elementary component  with more than two vertices is peripherally 2-colorable, and one component
has at least two finite faces,
then we can apply Algorithm 3.3 in \cite{C24+} on each component to generate   binary coding for all perfect matchings
of $G'$ inducing an isometric embedding of $R(G')$ into $Q_{d}$  as a finite distributive lattice,
where  $d=\mathrm{idim}(R(G'))$ is the sum of number of finites faces of all peripherally 2-colorable components of $G'$. 

Therefore, if $G'$ is a plane weakly elementary bipartite graph
 whose each elementary component  with more than two vertices is peripherally 2-colorable, and 
 one component has at least two finite faces.
 Then there are two distinct binary codings of length $d$ for all perfect matchings  of $G'$
 such that one binary coding induces an isometric embedding of $R(G')$ into $Q_{d}$  as a finite distributive lattice,
 and the other induces  an isometric embedding of $R(G')$ into $Q_{d}$  as a daisy cube,
 where   $d=\mathrm{idim}(R(G'))$.

\begin{corollary}\label{C:TwoCodings}
Let $G$ be a peripherally 2-colorable graph that is not an even cycle.
Assume that $G$ has  a $\mathrm{RFD}(G_1,G_2, \ldots ,G_n)$
 associated with a sequence of finite faces $s_1, s_2, \ldots, s_n$, where $n \ge 2$,
 and $s_i$ is a reducible face of $G_i$ for $2 \le i \le n$.
Let $\mathcal{L}_n^{FDL}$ (respectively, $\mathcal{L}_n^{DC}$) be a binary coding for all perfect matchings of $G$ 
which induces an isometric embedding of $R(G)$ into a hypercube $Q_n$  as a finite distributive lattice 
(respectively, as a daisy cube).
Let $M$ be an arbitrary  perfect matching of $G$. 

(i) Let $x_M$ be the binary string of $M$ in $\mathcal{L}_n^{FDL}$
such that each position $i$ of $x_M$ corresponding to a finite face $s_i$ of $G$  for $1 \le i \le n$.
If all exterior handles on $\partial s_i$ are proper $M$-alternating along the clockwise orientation of $\partial G$, 
then the $i$-th position of $x_M$  is $1$, and $0$ otherwise.

(ii) Let $y_M$ be the binary string of $M$ in $\mathcal{L}_n^{DC}$
such that  each position $i$ of $y_M$ corresponding to a finite face $s_i$ of $G$ for $1 \le i \le n$.
If $s_i$ is a finite face of $G$
such that $M$ does not contain end edges of any exterior handles on $\partial s_i$,  
then $i$-th position of $y_M$ is $0$, and $1$ otherwise.
 \end{corollary}	
	
\proof Let $G$ be a peripherally 2-colorable graph that is not an even cycle.
We observe that any facial cycle of $G$ is a union of odd length handles that
are interior and exterior  alternately along the clockwise orientation of the facial cycle 
such that any two consecutive ones  have exactly one vertex in common.

(i) Assume that $G$ is properly 2-colored such that any two adjacent vertices are assigned different colors.
Then either all exterior handles on the periphery of $s_i$  start from white vertices and end with black vertices,
or all exterior handles on the periphery of $s_i$  start from black vertices and end with white vertices
along the clockwise orientation of the periphery of $G$.

If $n=2$, then $R(G)$ is a path on three vertices, and it is easy to check that the conclusion is trivial.
Let $n>2$. We proceed by mathematical induction on $n$.
Assume that the conclusion is true for any peripherally 2-colorable graph that is not an even cycle
and with less than $n$ finite faces.
 Let $M$ be a perfect matching of $G$. Let $x_M$ be the binary string of $M$ in $\mathcal{L}_n^{FDL}$
such that each position $i$ of $x_M$ corresponds to a finite face $s_i$ of $G$  for $1 \le i \le n$.
By induction hypothesis,
for $1 \le i \le n-1$,  if all exterior handles on $\partial s_i$ are proper $M$-alternating along the clockwise orientation of $\partial G$, 
then the $i$-th position of $x_M$  is $1$, and $0$ otherwise.
By the proof of Theorem 3.4 \cite{C24+},  we can see that
if the exterior handle on $\partial s_n$ is proper $M$-alternating along the clockwise orientation of $\partial G$, 
then the $n$-th position of $x_M$  is $1$, and $0$ otherwise.

(ii)  Theorem \ref{T:ArbitrayFiniteFace} (i),  any perfect matching 
$M$ of $G$ either contains end edges of all exterior handles of $s_i$ or 
does not contain end edges of any exterior handles of $s_i$.
Let $y_M \in \mathcal{L}_n^{DC}$. Then by Theorem \ref{T:ArbitrayFiniteFace} (ii), if 
$M$ does not contain end edges of any exterior handles on $\partial s_i$, 
then the $i$-th position of $y_M$ is $0$,  and $1$ otherwise.
\qed\\	
	
By Corollary \ref{C:TwoCodings}, 
we can see that switching two color classes of a proper 2-coloring on the vertex set of $G$ results in swapping 0s and 1s for each binary string
in $\mathcal{L}_n^{FDL}$, but does not affect any binary string in $\mathcal{L}_n^{DC}$.


  

\end{document}